\documentclass[10pt,twoside]{article}
\usepackage{graphicx}
\usepackage{amsmath}
\usepackage{Latex-document}

\title{\bf  Topology of Singular Algebraic Varieties\vskip 6mm}

\author{B. Totaro\vspace*{-0.5cm}\thanks{Department of Mathematics and
Mathematical Statistics, University of Cambridge, Wilberforce
Road, Cambridge CB3 0WB, UK. E-mail: b.totaro@dpmms.cam.ac.uk}}

\date{\vspace{-8mm}}

\def\Z{\text{\bf Z}}
\def\Q{\text{\bf Q}}
\def\R{\text{\bf R}}
\def\C{\text{\bf C}}
\def\P{\text{\bf P}}
\def\F{\text{\bf F}}
\def\arrow{\rightarrow}

\def\imp{\Rightarrow}

\def\str{\text{str}}

\begin{document}

\maketitle

\newtheorem{theorem}{Theorem}[section]
\newtheorem{corollary}[theorem]{Corollary}
\newtheorem{lemma}[theorem]{Lemma}

\thispagestyle{first} \setcounter{page}{533}

\begin{abstract}

\vskip 3mm

I will discuss recent progress by many people in the program of
extending natural topological invariants from manifolds to
singular spaces. Intersection homology theory and mixed Hodge
theory are model examples of such invariants. The past 20 years
have seen a series of new invariants, partly inspired by string
theory, such as motivic integration and the elliptic genus of a
singular variety. These theories are not defined in a topological
way, but there are intriguing hints of their topological
significance.

\vskip 4.5mm


\end{abstract}

\vskip 12mm

\section{Introduction} \label{section 1}\setzero
\vskip-5mm \hspace{5mm }

The most useful fact about singular complex algebraic varieties is
Hironaka's theorem that there is always a resolution of
singularities \cite{Hironaka}.  It has long been clear that the
non-uniqueness of resolutions poses a difficulty in many
applications. Many different methods have been used to get around
this difficulty so as to define invariants of singular varieties.
One approach is to try to describe the relation between any two
resolutions, leading to ideas such as cubical hyperresolutions
\cite{HC} and the weak factorization theorem (\cite{AKMW},
\cite{Wlodarczyk}). Another idea, coming from minimal model
theory, is to insist on the special importance of crepant
resolutions, and more generally to emphasize the role of the
canonical bundle. Recently the interplay between these two
approaches has been very successful, as I will describe.

The recent methods tend to be more roundabout than the direct
topological definition of intersection homology groups. It is
tempting to try to define suitable generalizations of intersection
homology groups in order to ``explain'' various results below
(3.2, 3.4, 4.1, 5.2).

\section{The weight filtration} \label{section 2}
\setzero\vskip-5mm \hspace{5mm }

Deligne discovered a remarkable structure on the rational
cohomology of any complex algebraic variety, not necessarily
smooth or compact: the weight filtration \cite{Deligneweight}.
This filtration expresses the way in which the cohomology of any
variety is related to the cohomology of smooth compact varieties.
It is a deep fact that the resulting filtration is well-defined.
For example, an immediate consequence of the well-definedness of
the weight filtration on cohomology with compact support is the
following fact, originally conjectured by Serre (\cite{Durfee},
\cite{DK}, \cite{Fulton}, p.~92).

{\bf Theorem 2.1. }\it For any complex algebraic variety $X$, not
necessarily smooth or compact, one can define ``virtual Betti
numbers'' $a_iX\in \Z$ for $i\geq 0$ such that

(1) if $X$ is smooth and compact, then the numbers $a_iX$ are the
Betti numbers $b_iX=\dim_{\Q} H^i(X,\Q)$;

(2)  for any Zariski-closed subset $Y\subset X$,
$a_iX=a_iY+a_i(X-Y)$. \rm

Using resolution of singularities, it is clear that the numbers
$a_iX$ are uniquely characterized by these properties. What is
less clear is the existence of such numbers. It follows, for
example, that if two smooth compact varieties $X$ and $Y$ can be
written as finite disjoint unions of locally closed subsets,
$X=\coprod X_i$ and $Y=\coprod Y_i$, with isomorphisms $X_i\cong
Y_i$ for all $i$, then $X$ and $Y$ have the same Betti numbers.
This is a topological property of algebraic varieties which has no
obvious analogue in a purely topological context.

The existence of the weight filtration, and consequently of the
virtual Betti numbers $a_iX$, was originally suggested by
Grothendieck's approach to the Weil conjectures on counting
rational points on varieties over finite fields. Indeed, the
number of $\F_q$-points of a variety clearly has an additive
property analogous to property (2) above. One proof of the
existence of the weight filtration for complex varieties reduces
the problem to the full Weil conjecture for varieties over finite
fields, proved by Deligne \cite{DeligneWeil}. Around the same
time, Deligne gave a more direct proof of the existence of the
weight filtration for complex varieties, using Hodge theory
\cite{DeligneHodge}. This is a classic example of the philosophy
that the deepest properties of algebraic varieties can often be
proved using either number theory or analysis, while they have no
``purely geometric'' proof.

In 1995, however, Gillet and Soul\'e gave a new proof of the
existence of the weight filtration \cite{GS}. They used ``only''
resolution of singularities and algebraic $K$-theory, specifically
the Gersten resolution. As a result of their more geometric proof,
they were able to define the weight filtration on the integral
cohomology or $\F_l$-cohomology of a complex algebraic variety,
not only on rational cohomology.

To understand what this means, let me describe the weight
filtration for a smooth complex variety $U$, not necessarily
compact. Using resolution of singularities, we can write $U$ as
the complement of a divisor with normal crossings $D$ in some
smooth compact variety $X$. For $i\geq 0$, let $X^{(i)}$ be the
disjoint union of the $i$-fold intersections of divisors. Then
there is a spectral sequence
$$E_1=H^j(X^{(i)},k)\imp H^{i+j}_c(U,k)$$
for any coefficient ring $k$. The weight filtration on the
compactly supported cohomology of $U$ is defined as the filtration
associated to this spectral sequence. Gillet and Soul\'e show that
for any coefficient ring $k$, this filtration is an invariant of
$U$, independent of the choice of compactification $U$. This is
not at all clear from the known invariance of this filtration for
$k=\Q$.

In fact, Gillet and Soul\'e proved more: for any coefficient ring
$k$, the spectral sequence is an invariant of $U$ from the $E_2$
term on. For $k=\Q$, the spectral sequence degenerates at $E_2$,
but this is not true with coefficients in $\Z$ or $\F_l$. As a
result, for general coefficients $k$, the groups in the $E_2$ term
are interesting new invariants of $U$ which are not simply the
associated graded groups to the weight filtration. They satisfy
Mayer-Vietoris sequences, and so can be considered as a cohomology
theory on algebraic varieties.

I can now explain a new application of  the geometric proof that
the weight filtration is well-defined. Namely, one can try to
define the weight filtration not only for algebraic varieties. The
point is that resolution of singularities holds more generally,
for complex analytic spaces, and even for real analytic spaces.
Gillet and Soul\'e's construction of the weight filtration uses
algebraic $K$-theory as well as resolution of singularities, and
it is not clear how to adapt the argument to an analytic setting.
But Guillen and Navarro Aznar improved Gillet and Soul\'e's
argument so as to construct the weight filtration using only
resolution of singularities \cite{GNA}. The details of their
argument use their idea of ``cubical hyperresolutions'' \cite{HC}.

Using the method of Guillen and Navarro Aznar, I have been able to
define the weight filtration for complex and real analytic spaces.
In more detail, let us define a compactification of a complex
analytic space $X$ to be a compact complex analytic space
$\overline{X}$ containing $X$ as the complement of a closed
analytic subset. Of course, not every complex analytic space has a
compactification in this sense. We say that two compactifications
of $X$ are equivalent if there is a third which lies over both of
them.

{\bf Theorem 2.2. }\it Let $k$ be any commutative ring. Then the
compactly supported cohomology $H^*_c(X,k)$ has a well-defined
weight filtration for every complex analytic space $X$ with an
equivalence class of compactifications. \rm

Any algebraic variety comes with a natural equivalence class of
compactifications, but in the analytic setting this has to be
considered as an extra piece of structure. On the other hand, the
theorem says that the weight filtration is well-defined on all
compact complex analytic spaces, with no extra structure needed.

For real analytic spaces, one has the difficulty that there is no
natural orientation, unlike the complex analytic situation. This
is not a problem if one uses $\F_2$-coefficients, and therefore
one can prove:

{\bf Theorem 2.3. }\it For every real analytic space $X$ with an
equivalence class of compactifications, the compactly supported
cohomology of the space $X(\R)$ of real points with $\F_2$
coefficients has a well-defined weight filtration. \rm

In particular, one can define virtual Betti numbers $a_iX$ for a
real analytic space $X$ with an equivalence class of
compactifications, the integers $a_iX$ being the usual
$\F_2$-Betti numbers in the case of a closed real analytic
manifold.

{\bf Example. }Let $X$ be the compact real analytic space obtained
by identifying two copies at the circle at a point, and let $Y$ be
the compact real analytic space obtained by identifying two points
on a single circle (the figure eight). It is immediate to compute
that $a_0X=1$ and $a_1X=2$, whereas $a_0Y=0$ and $a_1=Y$. The
interesting point here is that the spaces $X(\R)$ and $Y(\R)$ of
real points are homeomorphic. Thus the numbers $a_i$ for a compact
real analytic space are not topological invariants of the space of
real points. In a similar vein, Steenbrink showed that the weight
filtration on the rational cohomology of complex algebraic
varieties is not a topological invariant, using 3-folds
\cite{Steenbrink}.

Nonetheless, it seems fair to say that extending the weight
filtration and the virtual Betti numbers to complex and real
analytic spaces helps to bring out more of the topological meaning
of these invariants of algebraic varieties. A real analytic space
has in some ways a weak structure; for example, the classification
of closed real analytic manifolds up to isomorphism is the same as
the classification of closed differentiable manifolds up to
diffeomorphism. From this point of view, it is surprising that
compactified real analytic spaces have the extra structure of the
weight filtration on their $\F_2$-cohomology. It seems natural to
ask for an $\F_2$-linear abelian category of ``mixed motives''
associated to compactified real analytic spaces $X$, such that the
$\F_2$-cohomology groups of $X$ with their weight filtration are
determined by the mixed motive of $X$. On Beilinson's conjectured
abelian category of mixed motives in algebraic geometry, see for
example Jannsen \cite{Jannsenbook}, 11.3, and \cite{JannsenMurre};
on various approximations to this category, see the triangulated
categories defined by Hanamura \cite{Hanamura}, Levine
\cite{Levine}, and Voevodsky \cite{VSF}, and the abelian category
defined by Nori.

It should be much easier to define mixed motives for real analytic
spaces than to do so for algebraic varieties. In particular, one
might speculate that the mixed motive of a real analytic space
should not involve much more information than the weight spectral
sequence converging to its $\F_2$-cohomology (starting at $E_2$),
perhaps considered together with an action of the Steenrod
algebra. In low dimensions, one could hope for precise
classifications of mixed motives along these lines.

\section{Stringy Betti numbers} \label{section 3}
\setzero\vskip-5mm \hspace{5mm }

The following result of Batyrev's \cite{Batyrev} is related to his
famous result that two birational Calabi-Yau manifolds have the
same Betti numbers. The proof uses Kontsevich's idea of motivic
integration \cite{Kontsevich}, as developed by Denef and Loeser
\cite{DL}. To be precise, Batyrev's statement involves Hodge
numbers, but I will only state what it gives about Betti numbers.

{\bf Theorem 3.1. }\it Let $Y$ be a complex projective variety
with log-terminal singularities. Then one can define the ``stringy
Poincar\'e function'' $p_{\str}(Y)$, which is a rational function,
such that for any crepant resolution of singularities $\pi:X\arrow
Y$, the stringy Poincar\'e function of $Y$ is the usual Poincar\'e
polynomial of $X$. \rm

We recall Reid's important definitions which are used here.
 First, let $Y$ be any normal complex
variety such that the canonical divisor $K_Y$ is $\Q$-Cartier. By
Hironaka, $Y$ has a resolution of singularities $\pi:X\arrow Y$
such that the exceptional divisors $E_i$, $i\in I$, are smooth
with normal crossings. The discrepancies $a_i$ of $E_i$ are
defined by
$$K_X=\pi^*K_Y+\sum a_iE_i.$$
The variety $Y$ is defined to have log-terminal singularities if
and only if $a_i>-1$ for all $i$. A resolution $X\arrow Y$ is said
to be crepant if $K_X=\pi^*K_Y$.

Batyrev defines the stringy Poincar\'e function of $Y$ by the
formula:
$$p_{\str}(Y)=\sum_{J\subset I}p(E_J^0)\prod_{j\in J}
\frac{q-1}{q^{a_j+1}-1}.$$ Here $E^0_J$ is the open stratum of
$E_J:=\cap_{j\in J}E_j$, and $p(E_J^0)$ denotes the virtual
Poincar\'e polynomial of $E_J^0$, written as a polynomial in
$q^{1/2}$. Thus $p_{\str}(Y)$ is a rational function in $q^{1/2}$
for $Y$ Gorenstein, and in $q^{1/n}$ for some $n$ in general.

Batyrev's proof that the stringy Poincar\'e function of $Y$ is
independent of the choice of resolution, using motivic
integration, rests on the additivity properties of the virtual
Poincar\'e polynomial. Using our extension of virtual Betti
numbers to complex analytic spaces, we find:

{\bf Theorem 3.2. }\it The stringy Poincar\'e function can be
defined as a rational function for any compactified complex
analytic space with log-terminal singularities. For any crepant
resolution $X\arrow Y$ with $Y$ compact, the stringy Poincar\'e
function of $Y$ is the usual Poincar\'e polynomial of $X$. \rm

Likewise for real analytic spaces:

{\bf Theorem 3.3. }\it An $\F_2$-analogue of the stringy
Poincar\'e function can be defined as a rational function for
compactified real analytic spaces with log-terminal singularities.
For any crepant resolution $X\arrow Y$ with $Y$ compact, the
stringy Poincar\'e function of $Y$ is the usual Poincar\'e
polynomial of the $\F_2$-cohomology of $X$. \rm

In particular, this answers part of Goresky and MacPherson's
Problem 7 in \cite{GM}:

{\bf Corollary 3.4. }\it Given a compact real algebraic variety
$Y$, the $\F_2$-Betti numbers of any two projective IH-small
resolutions of $Y$ are the same. \rm

This uses the relation between IH-small resolutions and crepant
resolutions, which I worked out in \cite{Totaro} using results of
Kawamata \cite{KMM} and Wisniewski \cite{Wisniewski}. In the
complex situation, the corollary (for Betti numbers with any
coefficients) has a more direct proof, since the Betti numbers of
any small resolution of $Y$ are equal to the dimensions of the
intersection homology groups of $Y$. It is not yet known whether
one can define a new version of intersection homology groups with
$\F_2$-coefficients which would be self-dual for all compact real
analytic spaces. A possible framework for defining such a theory
has been set up by Banagl \cite{Banagl}.

\section{The elliptic genus of a singular variety} \label{section 4}
\setzero\vskip-5mm \hspace{5mm }

I found that any characteristic number which can be extended from
smooth compact complex varieties to singular varieties, compatibly
with small resolutions, must be a specialization of the elliptic
genus \cite{Totaro}. It was then an important problem to define
the elliptic genus for singular varieties. This was solved in a
completely satisfying way by Borisov and Libgober \cite{BL}:

{\bf Theorem 4.1. }\it Let $Y$ be a projective variety with
log-terminal singularities. Then one can define the elliptic genus
of $Y$, $\varphi(Y)$, such that for any crepant resolution
$X\arrow Y$, we have $\varphi(Y)=\varphi(X)$. \rm

Here is Borisov and Libgober's definition of $\varphi(Y)$. Let
$\pi:X\arrow Y$ be a resolution whose exceptional divisors $E_k$
have simple normal crossings, and let $a_k$ be the discrepancies
as in section 3. Formally, let $y_l$ denote the Chern roots of $X$
so that $c(TX)=\prod_l(1+y_l)$, and let $e_k$ be the cohomology
classes on $X$ of the divisors $E_k$. Then $\varphi(Y)$ is the
analytic function of variables $z$ and $\tau$ defined by
$$\varphi(Y)=\int_Y (\prod_l \frac{(\frac{y_l}{2\pi i})\theta
(\frac{y_l}{2\pi i}-z)\theta'(0)}{\theta(-z)\theta(\frac{y_l}{2\pi
i})}) \times (\prod_k \frac{\theta(\frac{e_k}{2\pi
i}-(\alpha_k+1)z)\theta (-z)}{\theta(\frac{e_k}{2\pi
i}-z)\theta(-(\alpha_k+1)z)}),$$ where $\theta(z,\tau)$ is the
Jacobi theta function. The proof that $\varphi(Y)$ is independent
of the choice of resolution for log-terminal $Y$ uses the weak
factorization theorem of Abramovich, Karu, Matsuki, and Wlodarczyk
(\cite{AKMW}, \cite{Wlodarczyk}).

In the spirit of earlier sections, the singular elliptic genus
extends to compact complex analytic spaces with log-terminal
singularities. But it remains a mystery how to define the elliptic
genus for some topologically defined class of singular spaces that
would include singular analytic spaces with log-terminal
singularities.

\section{Possible characteristic numbers for real analytic spaces}
\label{section 5} \setzero\vskip-5mm \hspace{5mm }

In my paper \cite{Totaro}, in trying to define characteristic
numbers for singular complex varieties, it was very helpful to
require that these numbers are compatible with IH-small
resolutions, as Goresky and MacPherson had suggested (\cite{GM},
Problem 10). The problem thereby becomes more precise: it may be
possible to show that some characteristic numbers extend to
singular varieties and some do not. This can help to suggest
valuable invariants for singular varieties, such as Borisov and
Libgober's elliptic genus for singular varieties, even if one is
not a priori interested in IH-small resolutions. (The same
comments apply to crepant resolutions.)

With this in mind, we here begin to analyze which characteristic
numbers can be defined for real analytic spaces, or for
topological spaces with similar singularities, compatibly with
IH-small resolutions.
 In the complex situation, the fundamental example
of a singularity with two different IH-small resolutions is the
3-fold node; one says that the two IH-small resolutions are
related by the simplest type of ``flop.'' Likewise, in the real
situation, the real 3-fold node has two different IH-small
resolutions. For convenience, let us say that two closed
$n$-manifolds are related by a ``real flop'' if they are the two
different IH-small resolutions $X_1$ and $X_2$ of a singular space
with singular set of real codimension 3 that is locally isomorphic
to the product of the 3-fold node with an $(n-3)$-manifold.

Let us first consider characteristic numbers for unoriented
spaces. By Thom,
 the bordism ring $MO_*$ for unoriented manifolds is detected by
Stiefel-Whitney numbers. Therefore we can ask which
Stiefel-Whitney numbers (meaning $\F_2$-linear combinations of
Stiefel-Whitney monomials) are unchanged under real flops. Or,
more or less equivalently: what is the quotient of the bordism
ring $MO_*$ by the ideal of real flops $X_1-X_2$, for $X_1$ and
$X_2$ as above? There is a good answer:

{\bf Theorem 5.1. }\it The $\F_2$-vector space of Stiefel-Whitney
numbers which are invariant under real flops of $n$-manifolds is
spanned by the numbers $w_1^iw_{n-i}$ for $0\leq i\leq n$, or
equivalently by the numbers $w_1^{n-2i}v_i^2$ for $0\leq i\leq
n/2$, modulo those Stiefel-Whitney numbers which vanish for all
$n$-manifolds. Here $v_i=v_i(w_1,w_2,\ldots)$ denotes the Wu
class. The dimension of this space of invariant Stiefel-Whitney
numbers, modulo those which vanish for all $n$-manifolds, is 0 for
$n$ odd and $\lfloor n/2\rfloor +1$ for $n$ even. The quotient
ring of $MO_*$ by the ideal of real flops is isomorphic to:
$$\F_2[\R\P^2,\R\P^4,\R\P^8,\ldots]/((\R\P^{2^a})^2=(\R\P^2)^{2^a}
\text{ for all }a\geq 2).$$ \rm

This class of Stiefel-Whitney numbers has occurred before, in
Goresky and Pardon's calculation of the bordism ring of locally
orientable $\F_2$-Witt spaces \cite{GP}. To be precise, the latter
ring coincides with the above ring in even dimensions but is also
nonzero in odd dimensions. Goresky defined a Wu class $v_i$ in
intersection homology for $\F_2$-Witt spaces \cite{Goresky}, so
that the square $v_i^2$ lives in ordinary homology, and the
characteristic numbers for locally orientable $\F_2$-Witt spaces
$Y$ are obtained by multiplying these homology classes by powers
of the cohomology class $w_1$.

This does not explain the invariance of these Stiefel-Whitney
numbers for real flops, however. The problem is that the 3-fold
node is not an $\F_2$-Witt space. (Topologically, it is the cone
over $S^1\times S^1$, whereas the cone over an even-dimensional
manifold is a Witt space if and only if the homology in the middle
dimension is zero.) That is, the standard definition of
intersection homology is not self-dual on a space with 3-fold node
singularities. This again points to the problem of defining a new
version of intersection homology with $\F_2$ coefficients which is
self-dual on real analytic spaces. That should yield an $L$-class
in the $\F_2$-homology of such a space, which we can also identify
with the square of the Wu class, and which therefore should allow
the above characteristic numbers to be defined for a large class
of real analytic spaces. There are related results by Banagl
\cite{BanaglL}, for spaces which admit an extra ``Lagrangian''
structure.

We now ask the analogous question for oriented singular spaces:
what characteristic numbers can be defined, compatibly with
IH-small resolutions? We could begin by asking for the quotient
ring of the oriented bordism ring $MSO_*$ by oriented real flops
$X_1-X_2$, defined exactly as in the unoriented case ($X_1$ and
$X_2$ are the two small resolutions of a family of real 3-fold
nodes), except that we require $X_1$ and $X_2$ to be compatibly
oriented. It turns out that this is not enough: all Pontrjagin
numbers are invariant under oriented real flops, whereas they can
change under other changes from one IH-small resolution to
another, such as complex flops (between the two small resolutions
of a complex family of complex 3-fold nodes). By considering both
real and complex flops, we get a reasonable answer:

{\bf Theorem 5.2. }\it The quotient ring of $MSO_*$ by the ideal
generated by oriented real flops and complex flops is:
$$\Z[\delta,2\gamma,2\gamma^2,2\gamma^4,\ldots ],$$
where $\C\P^2$ maps to $\delta$ and $\C\P^4$ maps to
$2\gamma+\delta^2$. This quotient ring is exactly the image of
$MSO_*$ under the Ochanine elliptic genus (\cite{Landweber},
p.~63). \rm

This result suggests that it should be possible to define the
Ochanine genus for a large class of compact oriented real analytic
spaces, or even more general singular spaces.

\label{lastpage}


\begin{thebibliography}{99}
%\bibitem{1} M. F. Atiyah, The geometry of classical particles, {\it
%Asian J. of Maths.} (to appear).
%\bibitem{2} M. F. Atiyah \& R. Bott, The moment map and equivariant
%cohomology, {\it Topology}, 23 (1984), 1--28.
%\bibitem{3} V. Guillemin \& S. Sternberg, {\it Supersymmetry and
%equivariant de Rhaml theory}, Springer, 2000.

\bibitem{AKMW} D.~Abramovich, K.~Karu, K.~Matsuki, and J.~Wlodarczyk,
Torification and factorization of birational maps,
math.AG/9904135.

\bibitem{Banagl} M.~Banagl, Extending intersection homology type invariants
to non-Witt spaces, {\it Memoirs of the AMS, }to appear.

\bibitem{BanaglL} M.~Banagl, The $L$-class of non-Witt spaces, to appear.

\bibitem{Batyrev} V.~Batyrev, Stringy Hodge numbers of varieties
with Gorenstein canonical singularities, {\it Integrable systems
and algebraic geometry (Kobe/Kyoto, 1997)}, 1--32, World
Scientific, 1998.

\bibitem{BL} L.~Borisov and A.~Libgober, Elliptic genera of singular
varieties, {\it Duke Math.\ J., }to appear.

\bibitem{DK} V.~Danilov and A.~Khovanskii, Newton polyhedra and an
algorithm for computing Hodge-Deligne numbers, {\it Math.\ USSR
Izv.\ }{\bf 29 }(1987), 279--298.

\bibitem{DeligneHodge} P.~Deligne, Th\'eorie de Hodge I, II, III,
{\it Proc. ICM 1970}, v. 1, 425--430; {\it Publ.\ Math.\ IHES
}{\bf 40 }(1972), 5--57; {\bf 44 }(1974), 5--78.

\bibitem{DeligneWeil} P.~Deligne, La conjecture de Weil I,
{\it Publ.\ Math.\ IHES }{\bf 43 }(1974), 273--308.

\bibitem{Deligneweight} P.~Deligne, Poids dans la cohomologie des
vari\'et\'es alg\'ebriques, Actes ICM Vancouver 1974, I, 79--85.

\bibitem{DL} J.~Denef and F.~Loeser, Germs of arcs on singular algebraic
varieties and motivic integration, {\it Invent.\ Math.\ }{\bf 135
}(1999), 201--232.

\bibitem{Durfee} A.~Durfee, Algebraic varieties which are a disjoint
union of subvarieties, {\it Geometry and topology: manifolds,
varieties and knots, }99--102, Marcel Dekker, 1987.

\bibitem{Fulton} W.~Fulton, Introduction to toric varieties, Princeton, 1993.

\bibitem{GS} H.~Gillet and C.~Soul\'e, Descent, motives, and $K$-theory,
{\it J.\ reine angew.\ Math.\ }{\bf 478 }(1996), 127--176.

\bibitem{Goresky} M.~Goresky, Intersection homology operations,
{\it Comment.\ Math.\ Helv.\ }{\bf 59 }(1984), 485--505.

\bibitem{GM} M.~Goresky and R.~MacPherson, Problems and bibliography
on intersection homology, {\it Intersection homology}, ed.\
A.~Borel, Birkh\"auser, 1984, 221--233.

\bibitem{GP} M.~Goresky and W.~Pardon, Wu numbers of singular spaces,
{\it Topology }{\bf 28 }(1989), 325--367.

\bibitem{GNA} F.~Guill\'en and V.~Navarro Aznar, Un crit\`ere d'extension
d'un foncteur d\'efini sur les sch\'emas lisses, math.AG/9505008.

\bibitem{HC} F.~Guill\'en, V.~Navarro Aznar, P.~Pascual, and F.~Puerta,
{\it Hyperr\'esolutions cubiques et descente cohomologique},
Lecture Notes in Mathematics 1335, Springer, 1988.

\bibitem{Hanamura} M.~Hanamura, Homological and cohomological motives
of algebraic varieties, {\it Invent.\ Math.\ }{\bf 142 }(2000),
319--349.

\bibitem{Hironaka} H.~Hironaka, Resolution of singularities of an algebraic
variety over a field of characteristic zero, {\it Ann.\ Math.\
}{\bf 79 }(1964), 109--326.

\bibitem{Jannsenbook} U.~Jannsen, {\it Mixed motives and algebraic
$K$-theory}, LNM 1400, Springer, 1990.

\bibitem{JannsenMurre} U.~Jannsen, Motivic sheaves and filtrations
on Chow groups, {\it Motives (Seattle, 1991)}, 245--302, AMS,
1994.

\bibitem{KMM} Y.~Kawamata, K.~Matsuda, and K.~Matsuki, Introduction
to the minimal model program, {\it Algebraic Geometry (Sendai,
1985)}, ed.\ T.~Oda, 283--360, Kinokuniya-North Holland, 1987.

\bibitem{Kontsevich} M.~Kontsevich, lecture at Orsay, 7 December 1995.

\bibitem{Landweber} P.~Landweber, Elliptic cohomology and modular
forms, {\it Elliptic curves and modular forms in algebraic
topology, }55--68, LNM 1326, Springer, 1988.

\bibitem{Levine} M.~Levine, {\it Mixed motives,} AMS, 1998.

\bibitem{Steenbrink} J.~Steenbrink, Topological invariance of the weight
filtration, {\it Indag.\ Math.\ }{\bf 46 }(1984), 63--76.

\bibitem{Totaro}B.~Totaro, Chern numbers for singular varieties and
elliptic homology, {\it Ann.\ Math.\ }{\bf 151 }(2000), 757--791.

\bibitem{VSF} V.~Voevodsky, A.~Suslin, and E.~Friedlander, {\it Cycles,
transfers, and motivic homology theories, }Princeton, 2000.

\bibitem{Wisniewski} J.~Wisniewski, On contractions of extremal rays
of Fano manifolds, {\it J.\ reine angew.\ Math.\ }{\bf 417
}(1991), 141--157.

\bibitem{Wlodarczyk} J.~Wlodarczyk, Combinatorial structures on toroidal
varieties and a proof of the weak factorization theorem,
math.AG/9904076.
\end{thebibliography}
\end{document}